\documentclass{amsart}

\usepackage{amsmath, amsfonts, amssymb, amscd}

\newtheorem{theo}{Theorem}[section]
\newtheorem{lem}[theo]{Lemma}

\newtheorem{prop}[theo]{Proposition}

   \newtheorem{example}[theo]{Example}
   \newtheorem{definition}[theo]{Definition}

   \newenvironment{pf}{\noindent{\it Proof. }}{$\square$\par\medskip}

   \newcommand{\R}{{\mathbb R}}
   \newcommand{\C}{{\mathbb C}}

   \newcommand{\Aut}{\operatorname{Aut}}

\newcommand{\E}{\mathcal E}

   \renewcommand{\=}{\overset{\text{def}}{=}}

   \newfam\msbfam
   \font\tenmsb=msbm10 scaled \magstep1 \textfont\msbfam=\tenmsb
   \font\sevenmsb=msbm7 scaled \magstep1 \scriptfont\msbfam=\sevenmsb
   \font\fivemsb=msbm5 scaled \magstep1
   \scriptscriptfont\msbfam=\fivemsb

\newcommand{\SU}{\operatorname{SU}}

   \newcommand{\U}{{\mathcal U}}

   \newcommand{\V}{{\mathcal V}}

   \def\sideremark#1{\ifvmode\leavevmode\fi\vadjust{
   \vbox to0pt{\hbox to 0pt{\hskip\hsize\hskip1em
   \vbox{\hsize3cm\tiny\raggedright\pretolerance10000
   \noindent #1\hfill}\hss}\vbox to8pt{\vfil}\vss}}}

   \title[On the localization principle  for  pseudoellipsoids]
{On the localization principle for the  automorphisms of pseudoellipsoids}

\begin{document}

  \author{ Mario Landucci
and Andrea Spiro}




\begin{abstract} We show that Alexander's extendibility theorem  for a  local automorphism
of the unit ball  is valid also for  a  local automorphism   $f$ of a  pseudoellipsoid
$\E^n_{(p_1, \dots, p_{k})} \= \{ z \in \C^n :  \sum_{j= 1}^{n - k}|z_j|^2 + |z_{n-k+1}|^{2 p_1} + \dots + |z_n|^{2 p_{k}} < 1\ \}$, provided
 that $f$  is defined on a region $\U \subset \E^n_{(p)}$
 such that: i) $\partial \U \cap \partial \E^n_{(p)}$ contains an open set of
 strongly pseudoconvex points; ii)  $\U \cap \{\ z_i = 0\ \} \neq \emptyset$ for  any
 $n-k +1 \leq i \leq n$. By the counterexamples we exhibit,  such  hypotheses  can be considered  as optimal.\end{abstract}

 \maketitle

\section{Introduction}
For a given $k$-tuple of integers $p = (p_1, \dots, p_k)$, with each $p_\ell \geq 2$, let us denote by
$\E^n_{(p_1, \dots, p_{k})}$ (or,  more simply,  $\E_{(p)}^n$)   the  pseudoellipsoid in $ \C^n$ defined by
 $$\E^n_{(p_1, \dots, p_{k})} \= \{\ z \in \C^n\ :\
 \sum_{j= 1}^{n - k}|z_j|^2 + |z_{n-k+1}|^{2 p_1} + \dots + |z_n|^{2 p_{k}} < 1\ \}\ .$$
When $k = 0$, we  assume  $\E^n_{(p)}$  to be  the unit ball $B^n = \{\ z\in  \C^n\ : \ |z| < 1\ \}$.
Now, let us consider the following definition.\par
\medskip
\begin{definition}\label{localautomorphisms}{\rm
We call {\it local automorphism of $\E^n_{(p)}$\/} any
 biholomorphic map $f: \U_1 \subset \E^n_{(p)} \to \U_2 \subset \E^n_{(p)}$ between two connected open subsets of $\E^n_{(p)}$ such that:
 \begin{itemize}
 \item[a)] each of the intersections $\partial  \U_i \cap \partial  \E^n_{(p)}$, $i = 1,2$, contains a boundary  open set
   $\Gamma_i \subset \partial \E^n_{(p)}$;
  \item[b)]ÿ there exists at least one sequence $\{\ x_k\ \} \subset \U_1$ which converges to a point $x_o \in \Gamma_1$, which is  not a limit point of  $\partial \U_1 \cap  \E^n_{(p)}$, and so that $\{\ f(x_k)\ \}$
  converges  to a point $\hat x_o \in \Gamma_2$, which is  not a limit point of  $\partial \U_2 \cap  \E^n_{(p)}$.
 \end{itemize}
We say that {\it a  local automorphism  $f: \U_1 \subset \E^n_{(p)} \to \U_2 \subset \E^n_{(p)}$ extends  to a global automorphism of $\E^n_{(p)}$\/} if there exists
 some $F \in \Aut(\E^n_{(p)})$ such that $F|_{\U_1 \cap \E^n_{(p)}} = f|_{\U_1 \cap \E^n_{(p)}}$.}
 \end{definition}
 \medskip
 By a celebrated theorem of Alexander and  its generalization obtained  by Rudin (\cite{Al,Ru}),  when  $ \E^n_{(p)}=
 B^n$, any local automorphism extends to a global one. This crucial extendibility result
 is often quoted as {\it localization principle for the automorphisms of} $B^n$ and it
has been  extended or established under  different but similar  hypotheses,   for a wide class of  domains
besides the unit balls
(see e. g. \cite{DS, Pi, Pi1}). On the other hand, even if  it is known that the pseudoellipsoids $\E^n_{(p)}$
share many useful  properties with  $B^n$   for what concerns the global automorphisms and the proper holomorphic maps (see f.i. \cite{We, La, LS, DS}),
some simple examples show  that Alexander's  theorem
 cannot be true in full generality for  a pseudoellipsoid $\E^n_{(p)}$   different from $B^n$ (see e.g. Example \ref{examples}  below). \par
 \medskip
Nonetheless, for each $\E^n_{(p)}$,  it is possible to determine, precisely and in an efficient way,
the class of local automorphisms  that   can be extended  to  global ones. In this short note we give a characterization of such local automorphisms by means of the  following generalization of Alexander's  theorem. \par
\medskip
\begin{theo} \label{maintheorem} Let  $f: \U_1 \subset \E^n_{(p)}   \to  \U_2 \subset\E^n_{(p)}$ be a local   automorphism of a pseudoellipsoid $\E^n_{(p)}$,
  with $p = (p_1, \dots, p_k)$,  and  satisfying the following two conditions:
  \begin{itemize}
  \item[i)] there exists a sequence $\{x_i \}$ as in (b) of Definition \ref{localautomorphisms}, whose
  limit point $x_o \in \partial \E_{(p)}$  is Levi non-degenerate;
  \item[ii)]
  for any
  $ n-k+1 \leq i \leq n$, the intersection $\U_1\cap  \{\ z_i = 0\ \}$ is  not empty.
  \end{itemize}
  Then $f$ extends to a global automorphism  $f \in \Aut(\E^n_{(p)})$.
\end{theo}
\medskip
We point out that the set $\partial \E^n_{(p)} \cap \bigcup_{i = n-k+1}^n \{\ z_i = 0\ \}$ coincides with
the set of points of Levi degeneracy of $\partial \E^n_{(p)}$. So,  Theorem \ref{maintheorem} can be roughly stated
saying that   $f$ is globally extendible as soon as it admits an holomorphic
extension to some  open subset $\U \subset \E^n_{(p)}$,  which intersects   each of the hyperplanes containing the Levi degeneracy set
of $\partial \E^n_{(p)}$ and, at the same time,  the boundary $\partial \U$
contains an open set of strongly pseudoconvex points of $\partial \E^n_{(p)}$.\par
 From  next  Example \ref{examples},  it will be clear that such  hypotheses  can be considered  as optimal.\par
 \medskip
 The properties of the pseudoellipsoid  used in the proof  are basically just  two: (1) It admits a finite ramified covering over the unit ball; (2) Its automorphisms
 are  ``lifts"    of  the automorphisms of the unit ball that preserve the singular values of the covering. Since  (2)
 is a consequence of (1), it is reasonable  to expect  that  a similar  result
should be true for any arbitrary ramified covering of the unit ball. \par
About this more general problem,   we refer to \cite{KLS, KS} for  what concerns the classification  of the domains in $\C^2$ that admit a ramified holomorphic covering over $B^2$.\par
  \bigskip
 \section{On the automorphisms of the unit ball}
First of all, we need to recall some basic facts on the automorphisms of the unit ball. Let  us denote by $\hat \imath: \C^n \to \C P^n$ the canonical embedding
 $$\hat \imath: \C^n \to \C P^n\ ,\qquad \hat\imath(z) = \left[ \begin{matrix}
 z_1 \\   \vdots \\z_n \\ 1\end{matrix}\right]\ $$
 and let  ${\hat \C}^n = \hat \imath(\C^n) = \C P^n \setminus\{ [w]\ : w_{n+1} = 0\ \}$. We  recall that,
  via the embedding,   $B^n$
 corresponds  to the projective open set  $\hat B^n = \left\{\ [w]\in \C P^n\ :\ <w, w> < 0\ \right\}$
 where we denote by $<,>$ the pseudo-Hermitian inner product on $\C^{n+1}$ defined by
 $$<w, z> = \bar w^t \cdot I_{n,1} \cdot z\ ,\qquad \text{where}\ \ I_{n,1} \=
  \left( \begin{matrix} I_{n} & 0\\
 0 & -1 \end{matrix}\right) \eqno(2.1)$$
 It is also known that a holomorphic map $F: B^n \to B^n$ is an automorphism of $B^n$ if and only if
 the corresponding map $\hat F = \hat \imath \circ F \circ  \hat \imath^{-1}: \hat B^n \to \hat B^n$ is  a projective linear
 transformation which preserves the quadric $\partial \hat B^n = \{\ [w]\ :\ <w,w> = 0\ \}$  (see e.g. \cite{Ve}). This means
 that $\hat F$ is
 of the form
 $$\hat F([z]) = [{\mathbb A} \cdot z]\ ,\eqno(2.2)$$
 where $\mathbb A$ is a matrix  in $\SU_{n,1} $, i.e. such that  $\overline{\mathbb A}^t I_{n,1} {\mathbb A} =I_{n,1}$ and with $\det {\mathbb A} = 1$.\par
 The correspondence $F \mapsto \hat F = \hat \imath \circ F \circ  \hat \imath^{-1}$ gives an isomorphism
 between $\Aut(B_n)$ and $\SU_{n,1}/K$, where $ K = \left\{\ \ e^{i \frac{2\pi k}{n+1}} I_{n+1}\ , \ 0\leq k\leq n\ \ \right\}$.
 \par
 The identification of the elements of $\Aut(B^n)$  with the corresponding  projective linear transformations is often quite  useful,  for instance in order  to establish the following fact (see also \cite{We}, \S 6).\par
 \medskip
 \begin{lem}  \label{firstlemma} Let $F = (F_1, \dots, F_n) \in \Aut(B^n)$ be an automorphism such that
 $$F\left(B^n \cap \{\ z_i = 0\ \}\right) \subset \{\ z_i = 0\ \}\ \eqno(2.3)$$
 for all $n - k +1 \leq i \leq n$. Then the components $F_i$ are of the following form
$$F_j(z) = \frac{\sum_{\ell = 1}^{n-k}  A_j^\ell z_\ell + b_j}{\sum_{\ell = 1}^{n-k} c^\ell z_\ell + d}\ ,\qquad \text{for} \ 1 \leq j \leq n-k\ ,\eqno(2.4)$$
$$F_j(z) = e^{i \theta_j} z_j \frac{1}{\sum_{\ell = 1}^{n-k} c^\ell  z_\ell + d}\ ,\qquad \text{for} \ n-k+1 \leq j \leq n\ ,\eqno(2.5)$$
for some $\theta_j \in \R$  and where $A = (A^i_j)$,  $b = (b_j)$, $c = (c^\ell)$ and $d$ are so that  $ \left( \begin{matrix}  A & b\\
c  & d \end{matrix}\right) \in \SU_{n-k,1}$. In particular, the  maps $F_j$, $1 \leq j\leq n-k$,  coincide with
 the component of an element of $\Aut(B^{n-k})$,  while $\sum_{j = 1}^{n-k} c^j  z_j + d \neq 0$   for any $ z \in B^n$.
\end{lem}
\begin{pf} By hypothesis, the corresponding automorphism $\hat F  = \hat \imath \circ F \circ \hat \imath^{-1}\in \Aut(\hat B^n)$
maps all hyperplanes $H_i = \{\ [w] \in \C P^n\ : \ w_i = 0\ \}$ into themselves and hence
fixes their  poles  relative  to  the quadric $\partial \hat B^n$, i.e. fixes all the points
$$[e_i] = [0: \ldots : 0: \underset{i-th\ place}1: 0 : \ldots : 0]\ , \qquad n-k +1 \leq i \leq  n\ .$$
This implies that the matrix ${\mathbb A}$ which determines the projective transformation $\hat F$ is
of the form
$${\mathbb A}= \left( \begin{matrix}  A & \begin{matrix} 0 & \dots & 0 \end{matrix}  & b\\
 \begin{matrix} 0 \\ \vdots \\ 0 \end{matrix} & \smallmatrix e^{i \theta_{n-k+1}} &\phantom{\hdots} & 0 \\
 \phantom{\vdots} & \ddots &  \phantom{\vdots}  \\
  0 & \phantom{\hdots} &  e^{i \theta_{n}} \endsmallmatrix & \begin{matrix} 0 \\ \vdots \\ 0 \end{matrix}\\
c  &  \begin{matrix} 0 & \dots & 0 \end{matrix}   & d \end{matrix}\right)$$
where $A$, $b, c$  and $d \in \C$ are such that
${\mathbb A}' \=  \left( \begin{matrix}  A & b\\
c  & d \end{matrix}\right)$ belongs to $\SU_{n-k,1}$. From this, (2.4) and (2.5) follow immediately. The last claim follows from the fact that   the value  $\sum_{\ell = 1}^{n-k} c^\ell  z_\ell + d$  is the last homogeneous coordinate of the element  $[{\mathbb A}' \cdot (z_1: \ldots: z_{n-k}: 1] \in \C P^{n-k}$  and it is clearly different from $0$, since the map $[w] \mapsto [{\mathbb A}' \cdot w]$ is an automorphisms of $\hat B^{n-k} \subset \C P^{n-k} \setminus \{\ w_{n-k +1}\neq 0\  \}$.
 \end{pf}
 \bigskip
 \section{Proof of Theorem \ref{maintheorem}}
 First of all, we need to introduce the  following notation. For any $p = (p_1, \dots, p_k)$,   we will use the symbol $\pi^{(p)}$  to denote the map
 $$\pi^{(p)}: \C^n \to \C^n\ ,\qquad
 \pi^{(p)}(z) = (z_1, \dots, z_{n-k}, z_{n-k+1}^{p_1}, \dots, z_n^{p_{k}}) \ .$$
We recall that the restriction $\left.\pi^{(p)}\right|_{\E^n_{(p)}}$ gives a proper holomorphic map
 $\pi^{(p)}: \E^n_{(p)}  \longrightarrow B^n$.\par
 \medskip
Secondly, we need to recall a useful theorem by  Forstneric and Rosay (\cite{FR}).  Given a domain $D \subset \C^n$, we say that a boundary point $z_o\in\partial D$ {\it satisfies the condition $(P)$\/} if: 
\begin{itemize}
\item[--] $\partial D$ is of class $\mathcal C^{1 + \varepsilon}$ near  $z_o$ for some $\varepsilon >0$; 
\item[--] there exist a continuous negative plurisubharmonic function $\rho$ on $D$ and  a neighborhood $\U$ of $z_o$  so that  $\rho(z)\geq -c\ d(z,\partial D)$ at all points of $\U \cap D$ for some constant $c>0$.
\end{itemize}
Theorem 1.1  and some related remarks  of  \cite{FR}  can be summarized as follows.\par
\begin{theo} \label{FR} Let $h:D\rightarrow D'$ be a proper holomorphic map between  two domains of  $\C^n$ and  let $z_o\in\partial D$ be a point that satisfies the condition (P). \par
If   there exists a sequence  $\{z_j\} \subset D$  so that $\lim_{j \to \infty}Êz_j = z_o$ and  $ \lim_{j\to \infty} h(z_j) = \hat z_o$ for some   $\hat z_o \in\partial D'$  at which  $\partial D'$ is $\mathcal C^2$ and strictly pseudoconvex,   then $h$ extends continuously to all points of   neighborhood  $\V$ of $z_o$ in $\overline{D}$.
\end{theo}
 \medskip
 We may now prove the following lemma.\par
 \begin{lem} \label{secondlemma} Let $f: \U_1 \subset \E^n_{(p)}  \to \U_2 \subset \E^n_{(p)} $ be a local automorphism  of a pseudoellipsoid $\E^n_{(p)}$
  with $p = (p_1, \dots, p_k)$ and assume  that
  \begin{itemize}
  \item[i)] there exists a sequence $\{x_i \}$ as in (b) of Definition \ref{localautomorphisms}, whose
  limit point $x_o \in \partial \E_{(p)}$  is Levi non-degenerate;
  \item[ii)]
  for any
  $ n-k+1 \leq i \leq n$, the intersection $\U_1\cap  \{\ z_i = 0\ \}$ is  not empty.
  \end{itemize}
   Then, up to composition with a coordinate permutation
   $$(z_1, \dots, z_n) \mapsto (z_{\sigma(1)}, \dots, z_{\sigma(n)})\ ,\eqno(3.1)$$
    the map
  $f$ sends  the points of  the hyperplane $ \{\ z_i = 0\ \}$ into the same hyperplane for any
  $n - k + 1 \leq i \leq n$.
 \end{lem}
 \begin{pf} In all the following we will use the symbols $\Gamma_i$,  $x_o$ and $\hat x_o$ with the same
  meaning as in Definition  \ref{localautomorphisms}. \par
First of all, notice that  $\hat x_o \in \Gamma_2\subset\partial
\U_2$ satisfies the   condition (P)   and hence,
by Theorem \ref{FR}, for any
sufficiently small  ball $B_{\varepsilon}(\hat x_o)$,  centered at
$\hat x_o$ and of radius $\varepsilon$,  the holomorphic map
$f^{-1}: \U_2 \to \U_1$ extends continuously to all points of
$\overline{B_{\varepsilon}(\hat x_o)}  \cap \Gamma_2$. In
particular, we may  assume  that
$f^{-1}(\overline{B_{\varepsilon}(\hat x_o)  \cap \Gamma_2})$ is
contained in a neighborhood of $x_o = f^{-1}(\hat x_o)$  in
$\Gamma_1$ in which there are  no Levi degenerate point.\par Pick
a Levi non-degenerate point $\hat x_o' \in
\overline{B_{\varepsilon}(\hat x_o)}  \cap \Gamma_2$ and consider
a sequence $\{\hat x'_k\} \subset \overline{B_{\varepsilon}(\hat
x_o)}  \cap \U_2$ which converges to $\hat x'_o$. By construction,
the sequence  $\{ x_k' = f^{-1}(\hat x'_k)\} \subset \U_1$
converges to the Levi non-degenerate point $x'_o = f^{-1}(\hat
x'_o) \in \Gamma_1$. It follows  that, replacing $x_o$ by $x'_o$
and $\hat x_o$ by $\hat x_o'$ and by   Theorem \ref{FR}
applied to  $f$ and $f^{-1}$, there is no loss of generality if we
assume that $x_o$ and $\hat x_o$ are both Levi non-degenerate and
that, for any sufficiently small  $\varepsilon_1 > 0$,  the map
$f$ extends continuously to a map
$$f: \U_1\cup \left(\overline{B_{\varepsilon_1}(x_o)}  \cap \Gamma_1\right) \to  \U_2\cup \left(B_{\varepsilon}(\hat x_o) \cap  \Gamma_2\right)\ ,$$
which is an homeomorphism onto its image.
\par
\medskip
Since the complex Jacobian matrices   $\left.J\pi^{(p)}\right|_{x_o}$ and $\left.J\pi^{(p)}\right|_{\hat x_o}$  are of maximal rank   (recall that $x_o$ and $\hat x_o \in \partial \E^n_{(p)}$ are both  Levi non-degenerate), from  the fact that  $x_o$ is   not a limit point of $\partial \U_1 \cap \E^n_{(p)}$
 and by  the continuity of $f$ and $f^{-1}$ around $x_o$ and $\hat x_o$, respectively, we may choose
 $\varepsilon_1$ and $\varepsilon_2$ so that:
 \begin{itemize}
\item[a)] $\pi^{(p)}|_{B_{\varepsilon_1}( x_o)}$  and $\pi^{(p)}|_{B_{\varepsilon_2}( \hat x_o)}$  are both  biholomorphisms onto their  images;
\item[b)]  $\overline{f(B_{\varepsilon_1}(x_o)
ÿ\cap \U_1) }\subset B_{\varepsilon_2}(\hat x_o)$ and $f|_{B_{\varepsilon_1}(x_o)\cap \U_1}$ extends to an homeomorphism   between $\overline{B_{\varepsilon_1}(x_o)\cap \U_1}$ and $\overline{f(B_{\varepsilon_1}(x_o)\cap \U_1)}$ which induces an homeomorphism between $ B_{\varepsilon_1}(x_o) \cap \Gamma_1$ and $f(B_{\varepsilon_1}(x_o) \cap \Gamma_1) \subset \Gamma_2$;
\end{itemize}
Notice that, by definitions,  $x_o$ is not a limit point of $\partial \left(B_{\varepsilon_1}(x_o) \cap \U_1\right) \cap \E^n_{(p)}$ and,  by (b),  $\hat x_o$ is not a limit point of $\partial f \left(B_{\varepsilon_1}(x_o) \cap \U_1\right) \cap \E^n_{(p)}$.
So, if we set
$$\U'_1 \= B_{\varepsilon_1}(x_o)\cap \U_1\ , \quad \U'_2 \= f(\U'_1) \subset B_{\varepsilon_2}( \hat x_o)\ , \quad \V_i \= \pi^{(p)}(\U'_i) \ \ i = 1,2\ ,$$
 the maps
$$f|_{\U'_1}: \U'_1 \to \U'_2$$
and
 $$\tilde f = \left.\pi^{(p)} \circ f \circ  \pi^{(p)}{}^{-1}\right|_{\V_1}: \V_1 \subset B^n \longrightarrow \V_2 \subset B^n$$
 are  local automorphisms of  $\E^n_{(p)}$ and of the unit ball, respectively. \par
 By Rudin's
 generalization of Alexander's theorem (\cite{Ru}), this implies that $\tilde f$ extends to a global automorphism
 of $B^n$,  which we   denote by $\tilde f$ as well.
 By construction,  for any $z \in \U'_1 =  \pi^{(p)}{}^{-1}(\V_1)$, we have
 $$ \tilde f \circ \pi^{(p)}(z) = \pi^{(p)} \circ f(z)\ ,\eqno(3.2)$$
but since both sides have an holomorphic extension  on $\U_1$, we get that (3.2) must be  true also for any $z$ in such larger set.  \par
In particular,
$$J(\tilde f)|_{\pi^{(p)}(z)} \cdot J(\pi^{(p)})|_z =  J(\pi^{(p)})|_{f(z)} \cdot J(f)|_z \ ,\qquad
\text{for any} \ z \in \U_1\ .\eqno(3.3)$$
Since for any $z \in \U_1$, $\det J(f)|_z \neq 0$   and
$$\{\ J(\pi^{(p)})|_z = 0\ \} = \bigcup_{i = n - k +1}^n \{\ z_i = 0\ \}\ , \eqno(3.4)$$
equality  (3.3) implies that, for any  $n - k +1\leq i \leq n$ and $z \in  \U_1\cap  \{\ z_i = 0\ \}$,
the value of  $J(\pi^{(p)})|_{f(z)}$ is $0$.  By (3.4), this means that $ f\left(\U_1\cap  \{\ z_i = 0\ \}\right) $ is contained in the union $\bigcup_{j = n - k +1}^n \{\ z_j = 0\ \}$. Indeed, it is contained in exactly one of the hyperplanes $\{z_j = 0\}$, because $f$ is a biholomorphism and consequently  $ f\left(\U_1\cap  \{\ z_i = 0\ \}\right)$
is an irreducible analytic variety.  From this the conclusion follows.  \end{pf}
 We proceed  by defining a  rule that associates   an automorphism of
 $B^n$  with any local automorphism of a pseudoellipsoid  (see also \cite{We}, \S 6).
Given a local automorphism $f: \U \to \C^n$ of  $\E^n_{(p)} $, pick a point  $x_o \in \U \cap
\partial \E^n_{(p)}$ for which (b) of Definition \ref{localautomorphisms} holds
and determine a small ball $B_\varepsilon(x_o)$ centered in $x_o$ as in the proof of the previous
lemma. Then,  we denote by  $\tilde f  \in \Aut(B^n)$  the global automorphism of the unit ball
that extends  $\tilde f \= \pi^{(p)} \circ f \circ \pi^{(p)}{}^{-1}|_{\pi^{(p)}\left(\V\right)}$, with $\V \= B_\varepsilon (x_o) \cap \E^n_{(p)}$.  By the identity principle of the holomorphic maps, such automorphism $\tilde f$  depends only on  $f$  and  it will be called  {\it the (global) automorphism of $B^n$ associated with  $f$\/}.\par
  \medskip
  With the help of such correspondence, we may state the following  criterion  for extendibility of local automorphisms.\par
 \medskip
  \begin{prop} \label{mainproposition} A local   automorphism $f: \U_1 \subset  \E^n_{(p)} \to  \U_2 \subset  \E^n_{(p)}$ of a pseudoellipsoid $\E^n_{(p)}$, $p = (p_1, \dots, p_k)$,
 extends to a global automorphism $f \in \Aut(\E^n_{(p)})$ if and only if  its associated automorphism
 $\tilde f \in \Aut(B^n)$ satisfies (2.3) for any $ n-k+1 \leq i \leq n$,  up to composition with a
  permutation of those
 coordinates $z_{n-k + j}$, for which  the  integers $p_j$ are of the same value.
  \end{prop}
  \begin{pf} Assume that the local automorphism $f: \U \to \C^n$ extends to a global
  automorphism $f\in \Aut(\E^n_{(p)})$ and recall that,
  by construction,  the associated automorphism $\tilde f \in \Aut(B^n)$ satisfies (3.2)  at all points where $f$ is defined (in this case, at all points of $\E^n_{(p)}$).
  Then,  by Lemma \ref{secondlemma} and the fact that
  $\pi^{(p)}\left(\E^n_{(p)} \cap \{\ z_i = 0\ \}\right) = B^n \cap \{\ z_i = 0\ \}$,
   the equality  (3.3) implies that, up to a suitable permutation of coordinates,  $\tilde f$ satisfies (2.3) for any $ n-k+1 \leq i \leq n$.\par
  \medskip
  Conversely, assume that $f = (f_1, \dots, f_n) :\U_1 \subset  \E^n_{(p)} \to  \U_2 \subset  \E^n_{(p)}$  is  a
  local   automorphism  of $\E^n_{(p)}$ such that (up to a suitable permutation of coordinates)
  the associated automorphism  $\tilde f = (\tilde f_1, \dots, \tilde f_n)  \in \Aut(B^n))$ satisfies (2.3) for any $ n-k+1 \leq i \leq n$.
  From (2.4), (2.5) and (3.2), it follows that the component $f_{j}$ of $f$ are of the form
  $$f_j(z) = \frac{\sum_{\ell = 1}^{n-k} A_j^\ell z_\ell + b_j}{\sum_{\ell = 1}^{n-k} c^\ell z_\ell + d}\ ,\qquad \text{for} \ 1 \leq j \leq n-k\ ,\eqno(3.5)$$
$$f_{n-k+ j}(z) = e^{i \theta_j} z_j \frac{1}{\left(\sum_{\ell = 1}^{n-k} c^\ell  z_\ell + d\right)^{\frac{1}{p_j}}}\ ,\qquad \text{for} \ 1 \leq j \leq k\ ,\eqno(3.6)$$
for some fixed  definitions of the $p_j$-th roots  $w \mapsto w^{\frac{1}{p_j}}$. \par
>From (3.5) and (3.6) it follows immediately that $f$ coincides with a globally defined automorphism of
$\E^n_{(p)}$ (for the general expressions of the elements in $\Aut(\E^n_{(p)} )$ see \cite{We, La}).
  \end{pf}
  \bigskip
Now, Theorem \ref{maintheorem}  follows almost immediately. In fact,
 if  $f: \U_1 \subset  \E^n_{(p)} \to  \U_2 \subset  \E^n_{(p)}$  is a local   automorphism
  satisfying the hypothesis
 of the theorem, by Lemma \ref{secondlemma} and (3.2),  the associated
 automorphism $\tilde f\in \Aut(B^n)$ satisfies the hypothesis of Proposition \ref{mainproposition} and  the claim follows.\par
 \bigskip
 We conclude with the following simple  construction of non-extendible
 local automorphisms of pseudoellipsoids.\par
 \medskip
 \begin{example} \label{examples} {\rm Let $\tilde f \in \Aut(B^n)$ be an automorphism which does not satisfies (2.3) for some $n-k+1\leq j \leq n$.
 Pick a point $w_o \in \partial B \cap \{\ \prod_{j = n-k+1}^n z_j \neq  0 \ \}$
 so that also its image  $\tilde f(w_o)$ is in $\partial B \cap  \{\ \prod_{j = n-k+1}^n z_j \neq  0 \ \}$. Then, let $z_o \in \partial \E^n_{(p)}$
 so that  $\pi^{(p)}(z_o) = w_o$  and consider a connected neighborhood $\U$ of $z_o$ with the following two properties: a)
 $\pi^{(p)}|_{\U}$ is a biholomorphism between $\U$ and its image $\pi^{(p)}(\U)$; b) $\tilde f(\pi^{(p)}(\U))$   does not intersect  $  \left\{\ \prod_{j = n-k+1}^n z_j  =  0 \ \right\}$ (a sufficiently small neighborhood   $\U$ surely  satisfies both requirements). Then,  we may consider the map
 $$f: \U_1 = \U \cap \E^n_{(p)}  \to  \U_2 = f(\U) \cap  \E^n_{(p)} \ ,\qquad f \= \pi^{(p)}{}^{-1} \circ \tilde f \circ \pi^{(p)}\ . $$
By construction, $f$  is a local automorphism of $\E^n_{(p)}$ and its associated
 automorphism of $\Aut(B^n)$ is $\tilde f$.
 By the hypotheses  on $\tilde f$ and by Proposition \ref{mainproposition},
 $f$ cannot extend  to  a global automorphism of $\E^n_{(p)}$. }
 \end{example}

\bigskip
\bigskip
\font\smallsmc = cmcsc8
\font\smalltt = cmtt8
\font\smallit = cmti8
\hbox{\parindent=0pt\parskip=0pt
\vbox{\baselineskip 9.5 pt \hsize=3.1truein
\obeylines
{\smallsmc
Mario Landucci
Dip. Matematica Applicata  ``G. Sansone''
Universit\`a di Firenze
Via di Santa Marta 3
I-50139 Firenze
ITALY
}\medskip
{\smallit E-mail}\/: {\smalltt mario.landucci@unifi.it
}
}
\hskip 0.0truecm
\vbox{\baselineskip 9.5 pt \hsize=3.7truein
\obeylines
{\smallsmc
Andrea Spiro
Dip. Matematica e Informatica
Universit\`a di Camerino
Via Madonna delle Carceri
I-62032 Camerino (Macerata)
ITALY
}\medskip
{\smallit E-mail}\/: {\smalltt andrea.spiro@unicam.it}
}
}

    \end{document}